\RequirePackage{fix-cm}
\documentclass[twocolumn]{svjour3}          
\smartqed  

\newcommand{\be}{\begin{equation}}
\newcommand{\ee}{\end{equation}}

\usepackage{amsmath, amsfonts, amssymb, newlfont, graphicx}

\begin{document}

\title{Chaos in a non-autonomous nonlinear system describing asymmetric water wheels}


\titlerunning{Chaos in a non-autonomous nonlinear system describing asymmetric water wheels}        

\author{Ashish Bhatt and Robert A. Van Gorder}


\institute{A. Bhatt \at
Department of Mathematics, Universit{\"a}t Stuttgart \\
              \email{ashish.bhatt@mathematik.uni-stuttgart.de}    \\
R. A. Van Gorder \at
Mathematical Institute, University of Oxford, Andrew Wiles Building, Radcliffe Observatory Quarter, Woodstock Road, Oxford OX2 6GG United Kingdom \\
              \email{Robert.VanGorder@maths.ox.ac.uk}           
}

\date{Received: date / Accepted: date}

\maketitle

\begin{abstract}
We use physical principles to derive a water wheel model under the assumption of an asymmetric water wheel for which the water inflow rate is in general unsteady (modeled by an arbitrary function of time). Our model allow one to recover the asymmetric water wheel with steady flow rate, as well as the symmetric water wheel, as special cases. Under physically reasonable assumptions we then reduce the underlying model into a non-autonomous nonlinear system. In order to determine parameter regimes giving chaotic dynamics in this non-autonomous nonlinear system, we consider an application of competitive modes analysis. In order to apply this method to a non-autonomous system, we are required to generalize the competitive modes analysis so that it is applicable to non-autonomous systems. The non-autonomous nonlinear water wheel model is shown to satisfy competitive modes conditions for chaos in certain parameter regimes, and we employ the obtained parameter regimes to construct the chaotic attractors. As anticipated, the asymmetric unsteady water wheel exhibits more disorder than does the asymmetric steady water wheel, which in turn is less regular than the symmetric steady state water wheel. Our results suggest that chaos should be fairly ubiquitous in the asymmetric water wheel model with unsteady inflow of water.

\end{abstract}

\section{Introduction} \label{sec:intro}
\begin{figure}
\centering
\includegraphics[width =0.48\textwidth]{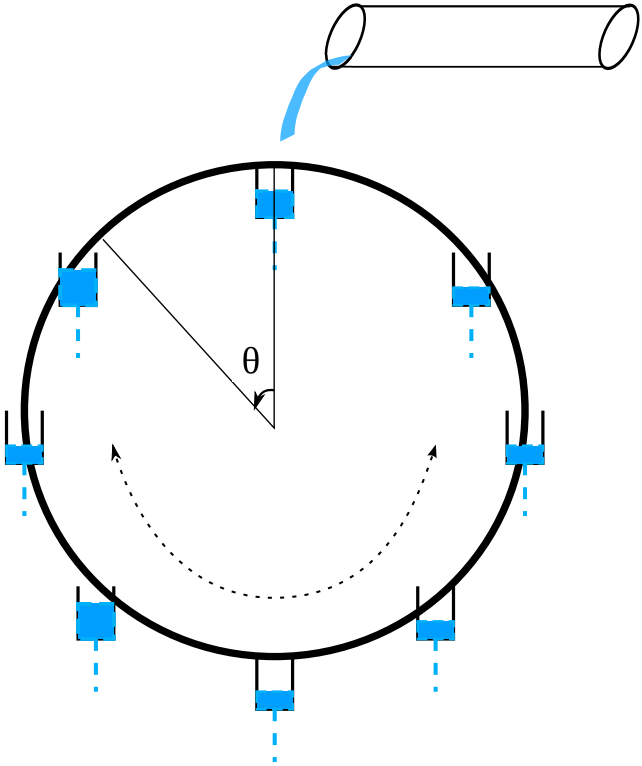}
\caption{A water wheel with water flowing into porous containers attached along the rim of the wheel.}
 \label{fig:waterwheel}
\end{figure}

A water wheel has several porous water containers attached along the rim of a 
wheel which rotates in a tilted plane. The angle $\theta \in [0,2\pi)$ is measured around the 
wheel in the counterclockwise direction. When the water is poured into the wheel 
at a certain angle, the containers start filling up and the wheel starts moving 
due to gravity. The water inflow and the loss of said water due to leakage create a 
seemingly random revolving motion where the wheel starts rotating in either 
direction at unpredictable instants of time; see Figure \ref{fig:waterwheel} for a schematic. 

To derive equations of motion of the water wheel, applying a standard conservation of mass argument \cite{Strogatz14} gives
\begin{align} \label{eq:mc}
\frac{\partial m}{\partial t} =q -Km -\omega \frac{\partial m}{\partial \theta},
\end{align}
while conservation of torque implies
\begin{align} \label{eq:Tc}
I\dot{\omega} = -\nu \omega +gR \int_0^{2\pi} m(\theta ,t)\sin \theta d\theta,
\end{align}
where dot denotes the usual time derivative, $\omega(t)$ is the angular velocity of the wheel, $m(\theta, t)$ is the mass distribution of the water around the rim of the wheel (defined in such a way that mass of the water between angles $\theta_1$ and $\theta_2$ is $M(t) = \int_{\theta_1}^{\theta_2} m(\theta,t) d\theta$), $q(\theta,t)$ is the water inflow (the rate at which water is pumped in above position $\theta$ at time $t$), $R$ is the radius of the wheel, $K$ is the leakage rate, $\nu$ is the rotational damping rate, $I$ is the moment of inertia of the wheel, and $g$ is the gravitational constant.

Since there is a little hope to find a solution of this 
coupled system in closed form, often one makes a transformation into a simpler yet more tractable problem which retains the salient features of the original problem. To this end, let us substitute the following Fourier series 
expansions of $m$ and $q$ in \eqref{eq:mc}-\eqref{eq:Tc}, we have
\begin{align} \label{eq:mqfe}
\begin{split}
m(\theta,t) = \sum_{n=0}^\infty [a_n(t) \sin n\theta +b_n(t) \cos n\theta],\\
q(\theta,t) = \sum_{n=0}^\infty [p_n(t) \sin n\theta +q_n(t) \cos n\theta].
\end{split}
\end{align}
When the water inflow $q$ is symmetric, no sine terms appear in the Fourier 
expansion of $q$. The presence of sine terms is a manifestation of 
\emph{asymmetric} inflow of water. As far as we are aware, all previous works have 
considered steady inflow of water \cite{Strogatz14,Lorenz63,Malkus72,KG92,MS06}, while we have assumed an unsteady inflow instead in the present work, for the sake of greater generality. This unsteady inflow results in the time-dependent 
harmonics $p_n$ and $q_n$ in the Fourier expansion of $q$.

Substituting \eqref{eq:mqfe} into \eqref{eq:mc}-\eqref{eq:Tc}, and using orthogonality of 
$\{\sin n\theta, \cos n\theta \}_{n=1}^\infty$ over the interval $[0,2\pi]$, we 
obtain the following system of coupled ODEs
\begin{align} \label{eq:aww-all-modes}
\begin{split}
\dot{a_n} &= n \omega b_n - Ka_n + p_n(t), \\
\dot{b_n} &= -n \omega a_n -Kb_n +q_n(t), \\
I\dot{\omega} &= -\nu \omega +\pi g R a_1,
\end{split}
\end{align}
for all $n=0,1,2,\ldots$. Retaining only the lowest modes, we obtain
\begin{align} \label{eq:aww-lowest-mode}
\begin{split}
\dot{a_1} &=  \omega b_1 - Ka_1 + p_1(t), \\
\dot{b_1} &= - \omega a_1 -Kb_1 +q_1(t), \\
I\dot{\omega} &= -\nu \omega +\pi g R a_1.
\end{split}
\end{align}
Making the change of variables $t =  \frac{\tau}{K}$, $\omega = Kx$, $a_1 = \frac{K\nu}{\pi g R}y$, $b_1 = \frac{-K \nu}{\pi g R}z +\frac{q_1}{K}$, $p_1 = \frac{K^2 \nu}{\pi g R}\mu$, $q_1 = \frac{K^2 \nu}{\pi g R}r$, $\nu = IK \sigma $ in \eqref{eq:aww-lowest-mode}, we obtain the following system
\begin{align} \label{eq:taww}
\begin{split}
\dot{x} &= \sigma (y-x), \\
\dot{y} &= r(\tau)x -y -xz +\mu(\tau), \\
\dot{z} &= xy -z +\dot{r}(\tau).
\end{split}
\end{align}
Here dot denotes a derivative with respect to the rescaled variable $\tau$. While the solution of system \eqref{eq:taww} does not provide a solution to the full system \eqref{eq:mc} - \eqref{eq:Tc}, it provides solution of the first harmonic 
in the Fourier expansion of $m(\theta,t)$ and hence is an approximation of the actual behavior of the full system.

Both symmetric and asymmetric water wheel models are known to possess chaos in 
certain parameter regimes. By eliminating all other possibilities, Lorenz showed 
in his seminal work \cite{Lorenz63} that the symmetric water wheel model has 
chaotic attractors for certain sets of system parameters. The asymmetric water wheel 
model has also been shown to be chaotic in certain parameter regimes in \cite{MS06}. 
Determining parameter values for which chaos occurs in a system has been mainly 
based on hunch and guesswork in the literature so far other than the most recent 
works based on the theory of competitive modes. This theory conjectures necessary conditions for a system to possess chaotic dynamics \cite{Pei06}. A very general 
quadratic system which encompasses steady symmetric and asymmetric water wheel models 
has been treated with competitive modes analysis in \cite{CG12} to reliably predict chaotic parameter 
regimes. In this paper, we use competitive modes analysis to predict chaotic 
parameter regimes in the system \eqref{eq:taww}. For a detailed overview and application of competitive modes analysis to a variety of ODE and PDE models, see \cite{LS06}. 

Many other seemingly disparate physical phenomena have been shown to have dynamical formulations which are slight variations of \eqref{eq:taww}. We mention some of them here without being exhaustive. It was shown in \cite{Malkus72} that the dynamics of thermal convection in a circular tube, with heat sources and sinks spread along its length, is equivalent to the Lorenz system. While thermal convection is a different physical phenomena than what we study in this paper, this example is similar in set-up and geometry as the heat sources and sinks are conceptually equivalent to the water inflow and leakage, and the circular tube analogous to the rim of the wheel. In \cite{Haken83}, the author has shown Maxwell-Block equations representing a laser model to be equivalent to the Lorenz system as well. In \cite{Robbins77}, the author has shown a model of disc dynamo to be equivalent to a special case of \eqref{eq:taww}. More recently, authors of \cite{CMT12} have shown a star shaped gas turbine to have a dynamical model similar to \eqref{eq:taww}. This turbine model is referred to as an augmented Lorenz model which consists of several Lorenz subsystems that share the angular velocity of the turbine as the central node.

The remainder of the paper is organized as follows. In Section \ref{sec:qa}, we analyze the dynamical systems \eqref{eq:aww-all-modes} and \eqref{eq:taww} qualitatively. Then we give an introduction to the theory of competitive modes and extend it to non-autonomous dynamical systems in Section \ref{sec:cm}. This theory is then applied to the system \eqref{eq:taww} to estimate chaotic parameter regimes in Section \ref{sec:cm2taww}. We discuss interesting findings in Section \ref{sec:conc}.

\section{Qualitative analysis of the water wheel dynamical system} \label{sec:qa}
Equation \eqref{eq:taww} represents the behavior of only the lowest modes. It is perhaps in order to discuss the behavior of higher modes in \eqref{eq:aww-all-modes}. It is rather straightforward to observe that
\begin{align} \label{eq:vol-cont}
\frac{\partial}{\partial a_n} \dot{a_n} + \frac{\partial}{\partial b_n} \dot{b_n} =-2K <0.
\end{align}
Therefore all volume elements in the phase space $(a_n,b_n)$ contract uniformly. Also from \eqref{eq:aww-all-modes} we get
$$\begin{aligned}
& -\frac{1}{2K} \frac{\partial}{\partial t} \left(a_n^2 +b_n^2\right) \\
& ~ = \left(a_n -\frac{p_n}{2K}\right)^2 +\left(b_n -\frac{q_n}{2K}\right)^2 -\left( \frac{p_n}{2K}\right)^2 -\left( \frac{q_n}{2K}\right)^2.
\end{aligned}$$
Since $K>0$, all trajectories in $(a_n,b_n)$ phase space get attracted toward the boundary of the evolving circle
\begin{align} \label{eq:evol-circle}
\left(a_n -\frac{p_n}{2K}\right)^2 +\left(b_n -\frac{q_n}{2K}\right)^2 =\left( \frac{p_n}{2K}\right)^2 +\left( \frac{q_n}{2K}\right)^2
\end{align}
centered at $\left(\frac{p_n}{2K},\frac{q_n}{2K} \right)$ with radius $\sqrt{\left(\frac{p_n}{2K}\right)^2 +\left(\frac{q_n}{2K}\right)^2}$.
Therefore, all the trajectories remain bounded so long as the modes $q_n$ and $p_n$ remain bounded and get attracted toward the evolving circle whose perimeter intersects the origin in $(a_n,b_n)$ phase space. 

It is sufficient to characterize the evolving circle to characterize the higher modes $a_n$ and $b_n$. To this end, let us rewrite the system \eqref{eq:aww-all-modes} as a system of integral equations
\begin{gather}
\begin{aligned} \label{eq:int_form}
a(t) &= e^{-Kt} \left( a(0) + \int_0^t e^{Ks}(n\omega(s) b(s) + p(s))~ds \right),\\
b(t) &= e^{-Kt} \left( b(0) + \int_0^t e^{Ks}(-n\omega(s) a(s) + q(s))~ds \right),
\end{aligned}
\end{gather}
where we have suppressed subscript $n$ to simplify the notation. Assuming $a^{(0)} = b^{(0)} =0$ and solving the integral equation \eqref{eq:int_form} iteratively, we find
\begin{gather} \label{eq:int_form_sol}
\begin{aligned}
a(t) &= e^{-Kt} [ a(0) \cos(n\Omega(t)) +b(0)\sin(n \Omega(t)) \\
&+\sum_{k=0}^\infty (-1)^k \left(n^{2k} P_{2k}(t) +n^{2k+1} Q_{2k+1}(t)\right)],\\
b(t) &= e^{-Kt} [ -a(0) \sin(n\Omega(t)) +b(0)\cos(n \Omega(t)) \\
&-\sum_{k=0}^\infty (-1)^k \left(n^{2k+1} P_{2k+1}(t) -n^{2k} Q_{2k}(t)\right)],
\end{aligned}
\end{gather}
where
$$\Omega(t) = \int_0^t \omega(s)~ds, 
$$
and
$$P_k(t) = \int_0^t \omega(s) P_{k-1}(s)~ds, 
$$
$$Q_k(t) = \int_0^t \omega(s) Q_{k-1}(s)~ds$$
for $k = 1,2,3,...$, while for $k=0$ we have
$$P_0(t) = \int_0^t e^{Ks} p(s)~ds,
$$
$$Q_0(t) = \int_0^t e^{Ks} q(s)~ds.$$

Clearly, the domain of convergence of the series occurring in \eqref{eq:int_form_sol} depends on the the angular frequency $\omega$. Solution \eqref{eq:int_form_sol} is suggestive of the following solution:
\begin{gather} \label{eq:int_form_sol2}
\begin{aligned}
a(t) &= e^{-Kt} [ \alpha \cos(n\Omega(t)) +\beta \sin(n \Omega(t))] +\int_0^t p(s)ds,\\
b(t) &= e^{-Kt} [ \beta \cos(n \Omega(t)) -\alpha \sin(n\Omega(t))] +\int_0^t q(s)ds.
\end{aligned}
\end{gather}
While there can be a more general solution, see \cite{KG92}, this solution form will serve our purposes. Therefore all the solutions converge toward a single attractor as $t \to \infty$. This attractor must lie inside the circle mentioned above and is determined by the character of $\omega$.

If we assume a steady (i.e., $q(\theta,t) = q(\theta)$) but asymmetric inflow, we obtain a \emph{steady asymmetric} water wheel model. In contrast to previously studied steady water wheel models, here we consider an unsteady water wheel model: system \eqref{eq:taww} represents a water wheel where water inflow is asymmetric and unsteady. As a consequence of this latter attribute of the water inflow, the ODE system \eqref{eq:taww} is non-autonomous. If we assume a symmetric inflow (no sine terms in the Fourier series of $q(\theta,t)$, \eqref{eq:mqfe}, and $\mu =0$ in \eqref{eq:taww}), we obtain an \emph{unsteady symmetric} water wheel model, instead.

Asymmetric inflow results in an asymmetric system: negating all the system 
variables $(x,y,z)$ in \eqref{eq:taww} results in a different system. While 
Malkus' symmetric water wheel is equivalent to the Lorenz system 
\cite{Strogatz14,Lorenz63}, \eqref{eq:taww} cannot be reduced to the Lorenz 
system \cite{MS06}. Nonetheless, like the Lorenz system, the dynamical system 
\eqref{eq:taww} contracts volume:
$$\nabla \cdot F = -\sigma -2 <0,$$
where $F = [\dot x, \dot y, \dot z]^T$ is the vector field associated with the 
dynamical system. The volume contraction implies that any solution of the system converges to an absorbing set in phase space. This, in turn, implies that the system has no quasi-periodic solutions 
and no repelling fixed points or closed orbits. Therefore all fixed points must 
be sinks or saddles.

\section{General competitive modes analysis} \label{sec:cm}
Because of its proximity and similar characteristics to the Lorenz system, we expect chaos in certain parameter regimes for the system \eqref{eq:taww}. Traditionally, one would search the parameter space $(\sigma, r(\tau), \mu(\tau))$ for a set of parameters which result in a chaotic attractor. Competitive modes analysis makes this process of searching for chaotic parameter regimes more methodical and organized. The theory helps delineate chaotic regimes in nonlinear autonomous dynamical systems, and can be used to predict chaos in the parameter regimes where mode frequencies are competitive.

As suggested in \cite{Pei06}, one writes a first 
order dynamical system as a second order system by differentiating with respect to 
time and then putting it into the form of an oscillator system. Given a general 
nonlinear autonomous system
\begin{align} \label{eq:nads}
\dot x_i = f_i(x_1,x_2, \ldots, x_n), \quad i =1,2,\ldots,n,
\end{align}
where $f_i \in C^1(\mathbb{R})$ and $\left|\frac{\partial f_i}{\partial x_j} \right|$ 
is bounded for all $j$, one can get the following second order system of 
differential equations by differentiating \eqref{eq:nads} with respect to $t$:
\begin{equation}\begin{aligned} \label{eq:dot-nads}
\ddot x_i &= \sum_{j=1}^n f_j \frac{\partial f_i}{\partial x_j}  \\
&= -x_i g_i(x_1,x_2, \ldots, x_n) \\
& \qquad + h_i(x_1,x_2, \ldots, x_{i-1}, x_{i+1}, 
\ldots, x_n),
\end{aligned}\end{equation}
for $i =1,2,\ldots,n$. Equation \eqref{eq:dot-nads} is analogous to a series of oscillators with their frequencies given by $g_i$, $i=1,2,\ldots,n$. The following conjecture gives necessary conditions for chaos in \eqref{eq:nads} as restrictions on these frequencies and functions.

\begin{conjecture} \label{conj:chaos-cond}
The conditions for a dynamical system to be chaotic are given below:
\begin{enumerate}
\item there exist at least two $g$'s in the system;
\item at least two $g$'s are competitive or nearly competitive, that is, there 
are $g_i \backsimeq g_j > 0$ at some $t$;
\item at least one of $g$'s is the function of evolution variables such as $t$; 
and
\item at least one of $h$'s is the function of the system variables.
\end{enumerate}
\end{conjecture}
The $x_i$'s for which all four conditions of the conjecture are satisfied are 
called \emph{competitive modes}.

However the system of interest, given by equation \eqref{eq:taww}, is 
non-autonomous. Such systems can be converted into autonomous systems by
enlarging the set of system variables. To this end, given a non-autonomous 
dynamical system
\begin{align*}
\dot x_i = f_i(x_1,x_2, \ldots, x_n,t), \quad i =1,2,\ldots,n,
\end{align*}
let us define $x_{n+1} = t$ to get the following autonomous system
\begin{align} \label{eq:nnds}
\dot x_i = f_i(x_1,x_2, \ldots, x_n,x_{n+1}), \quad i =1,2,\ldots,n+1.
\end{align}
Here again $f_i \in C^1(\mathbb{R})$, $\left|\frac{\partial f_i}{\partial x_j} \right|$ 
is bounded for all $j$, and $f_{n+1} = 1$. With these definitions we get from 
\eqref{eq:nnds} that the necessary second-order oscillator system corresponding to a non-autonomous first order system is given by
\begin{align*}
\ddot x_i &= \sum_{j=1}^{n+1} f_j \frac{\partial f_i}{\partial x_j}, \text{ for 
} i = 1,2,\ldots,n,\\
\ddot x_{n+1} &= \sum_{j=1}^{n+1} f_j \frac{\partial f_{n+1}}{\partial x_j} =0,
\end{align*}
because $\frac{\partial f_{n+1}}{\partial x_j} =0$ for all $j$. This system is 
equivalent to \eqref{eq:dot-nads} with one additional system variable. Therefore, making use of this modification,
one can use the theory of competitive modes to analyze non-autonomous systems. 

Competitive modes analysis has been used to predict chaotic regimes and custom design chaotic systems. Competitive modes were shown to result in a chaotic solution on a slow time scale in parametrically driven surface waves \cite{CG85}. In \cite{YYED06}, a competitive modes analysis was used to predict chaos in nonlinear dynamical systems and was employed to construct custom designed chaotic systems. Competitive modes analysis has been used to obtain parameter regimes giving chaos in a variety of nonlinear systems, such as the T system \cite{CG10}, a generalized Lorenz system \cite{VanGorder11}, a general chaotic bilinear system of Lorenz type \cite{CM7}, and the canonical form of the blue sky catastrophe \cite{bsc}, to name a few application. In \cite{CG12}, a competitive modes analysis was used to find general chaotic regimes for a quadratic dynamical system by relaxing some conditions of Conjecture \ref{conj:chaos-cond}.

\section{Competitive modes analysis of the water wheel system} \label{sec:cm2taww}
The water wheel model can be characterized by the solution behavior of \eqref{eq:aww-all-modes}. The qualitative behavior of higher modes was discussed in Section \ref{sec:qa} and the lowest modes will be treated with competitive modes analysis introduced in the previous section. We compare the three water wheel models (viz., nonsteady asymmetric, nonsteady symmetric, and steady asymmetric), by solving them numerically using the following initial conditions and parameter values (unless mentioned otherwise): 
\be \sigma = 5, \label{sigma}\ee
\be \mu(\tau) = 1, \label{mu}\ee
\be\begin{aligned}
r(\tau) & = 50 + 0.5 \sin(10 \tau) - 21.5 \tanh(20 - \tau)  \\
& \qquad - 28 \tanh(40 -\tau), \label{r}
\end{aligned}\ee
\be p_2(\tau) = \frac{\mu(\tau)}{100},\ee
\be q_2(\tau) = \frac{r(\tau)}{100},\ee
\be \{y(0),z(0),a_2(0),b_2(0)\} = \{0,r(0),0,0\} .\ee
For the steady water wheel, we instead set 
\be r(\tau) = r(30). \ee
We shall keep $x(0)$ as a free parameter.

Let us rationalize our choice for the function $r(\tau)$. First, let us remember that the function $r(\tau)$ is a scalar multiple of the first harmonic $q_1$ in the 
Fourier expansion of the water inflow rate $q(\theta,t)$. A smooth approximation of step function for 
$q_1$ corresponds to changing the water inflow rate abruptly. Secondly, we make the inflow more realistic by including small sinusoidal oscillations to model small fluctuations in the inflow. The function $r(\tau)$ thus represents small fluctuations and large changes in the water inflow rate. Similar arguments can be used to justify the choice of $p_2(\tau)$ and $q_2(\tau)$.

To employ the competitive modes analysis, differentiating system \eqref{eq:taww} with respect to time and comparing with 
\eqref{eq:dot-nads}, we obtain
\begin{align}
\begin{split}
g_1 &= -\sigma^2 - \sigma (r(\tau)  - z), \\
g_2 &= -1 -\sigma (r(\tau)-z)  +x^2, \\
g_3 &= -1 +x^2, \\
g_4 &= 0, \\
h_1 &= \sigma \mu(\tau) -(\sigma +\sigma^2)y, \\
h_2 &= -\mu(\tau) -\dot{\mu}(\tau) -(1+\sigma)r(\tau)x +(2+\sigma)xz, \\
h_3 &= -\dot{r}(\tau) +\ddot{r}(\tau) +\mu(\tau)x -(2 +\sigma)xy +r(\tau)x^2 +\sigma y^2, \\
h_4 &= 0.
\end{split}
\end{align}
Clearly, conditions 1, 3, and 4 of Conjecture \ref{conj:chaos-cond} are satisfied. In 
order for the $g$ functions to satisfy condition $2$ of the conjecture, solutions of 
the dynamical system \eqref{eq:taww} must satisfy at least one of the following conditions:
\begin{equation} \label{eq:mfs1}
g_1 - g_4 = 0 \iff -\sigma^2 - \sigma (r(\tau)  - z)=0, 
\end{equation}
\begin{equation} \label{eq:mfs2}
g_2 - g_4 = 0 \iff -1 -\sigma (r(\tau)-z)  +x^2=0,
\end{equation}
\begin{equation} \label{eq:mfs3}
g_3 - g_4 = 0 \iff -1 +x^2=0,
\end{equation}
\begin{equation} \label{eq:mfs4}
g_1 - g_2 =0 \iff 1 - \sigma^2 -x^2 =0,
\end{equation}
\begin{equation} \label{eq:mfs5}
g_2 - g_3 =0 \iff r(\tau) -z=0,
\end{equation}
\begin{equation} \label{eq:mfs6}
g_1 - g_3 =0 \iff 1 - \sigma^2 -x^2 -\sigma(r(\tau) -z)=0.
\end{equation}

In the phase space of system variables, equations of system \eqref{eq:mfs1}-\eqref{eq:mfs6} represent 
manifolds such as parabolas and lines for fixed parameter values. Notice that the equation $g_1-g_2=0$ does not represent a manifold in $\mathbb{R}^2$ if $\sigma >1$. Also, the distance between both, the two parabolas and the two lines, is exactly equal to $\sigma$. Moreover the parabolas intersect with the lines at $x=\pm 1$, which is the equation of the manifold $g_3-g_4=0$. Another useful observation is that the parameters $r(\tau)$ and $\mu(\tau)$ are scalar multiples of the symmetric and asymmetric modes, $q_1(\tau)$ and $p_1(\tau)$, respectively. Between these two parameters, only the former enters into \eqref{eq:mfs1}-\eqref{eq:mfs6}, and hence only the symmetric part of the inflow, $q_1(\tau)$, influences chaos as far as the competitive modes analysis is concerned. It is also useful to note that a non-autonomous system of dimension $n$ will in general have $n$ additional manifolds compared to an autonomous system of the same dimension. Therefore, the non-autonomous system has more $g$'s to meet the conditions of Conjecture \ref{conj:chaos-cond} than the autonomous system, potentially enlarging the chaotic parameter regime.

\begin{figure}
\centering
\includegraphics[width =0.45\textwidth]{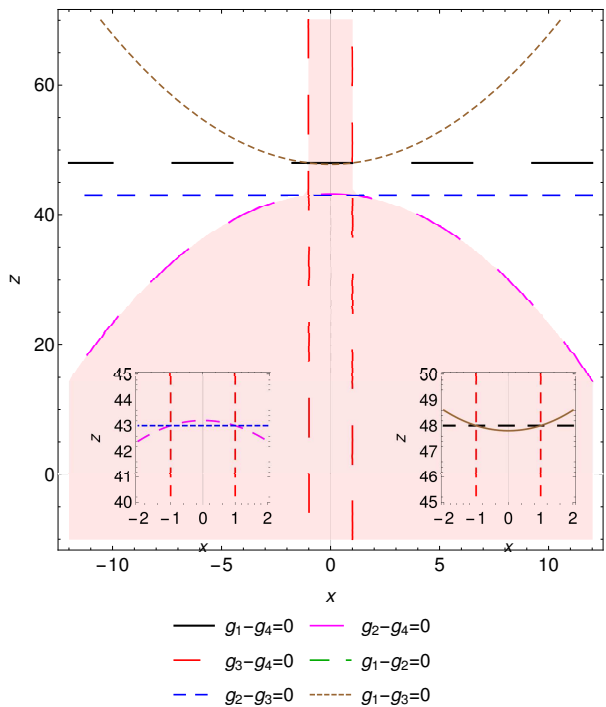}
\caption{Plots of the manifolds in the $x-z$ plane given by \eqref{eq:mfs1}-\eqref{eq:mfs6} at time $\tau=5$. The shaded area is where no competitive $(g_i,g_j)$ pair is positive in the plane.}
 \label{fig:gijpos}
\end{figure}

\begin{figure}
\centering
\includegraphics[width =0.35\textwidth]{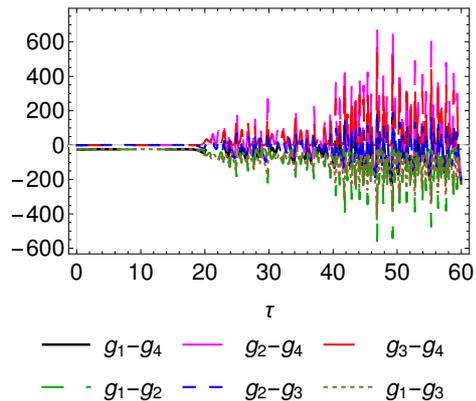}
\caption{Time series plots of the expressions $g_i -g_j$ in \eqref{eq:mfs1}-\eqref{eq:mfs6} evaluated along the numerical solution of the system \eqref{eq:taww}.}
 \label{fig:gi-gj}
\end{figure}

Figure \ref{fig:gijpos} shows a typical representation of the
manifolds listed in \eqref{eq:mfs1}-\eqref{eq:mfs6}, while Figure \ref{fig:gi-gj} gives corresponding time series plots. Figure \ref{fig:gijpos} is a snapshot of the evolving parameter-dependent manifolds. The maximum amplitude of the expressions $g_i - g_j$ increases with increasing value of the parameter $r(\tau)$ in Figure \ref{fig:gi-gj}.

Conjecture \ref{conj:chaos-cond} requires that the competitive $g$'s be positive. In order 
to verify this, we plot the $x-z$ region where none of the competitive 
$(g_i,g_j),i\neq j$ pair is positive with the same parameters as in the previous 
plot. This region is shown in Figure \ref{fig:gijpos}. Clearly, there are regions of 
the phase space where competitive $g$'s are positive for the chosen parameters. When the solution 
$(x,y,z)$ of the system \eqref{eq:taww} resides in these regions, positivity 
condition of Conjecture \ref{conj:chaos-cond} is satisfied. Clearly, manifolds $g_2=0$ and $g_3=0$ delineate the positive competitive region of the space from non-competitive region.

It is necessary for the solution $(x(\tau),y(\tau),z(\tau))$ to lie 
on at least one of the manifolds of Figure \ref{fig:gijpos} intermittently for the system 
to be chaotic. Solving \eqref{eq:taww} numerically (taking $x(0)=1$) and evaluating $g_i -g_j$, for all $i \neq 
j$, along the numerical solution, we obtain Figure \ref{fig:gi-gj}. For the chosen 
parameter values, the pairs $\{g_2,g_4\}, \{g_3,g_4\}$, and $\{g_2,g_3\}$ are 
competitive according to the time series plotted. No pair is competitive around the time 
interval $\tau \in [0,20]$, as expected, because $r(\tau) \approx 0$ in that 
interval.

\begin{figure*}
\centering
\includegraphics[width=0.9\textwidth]{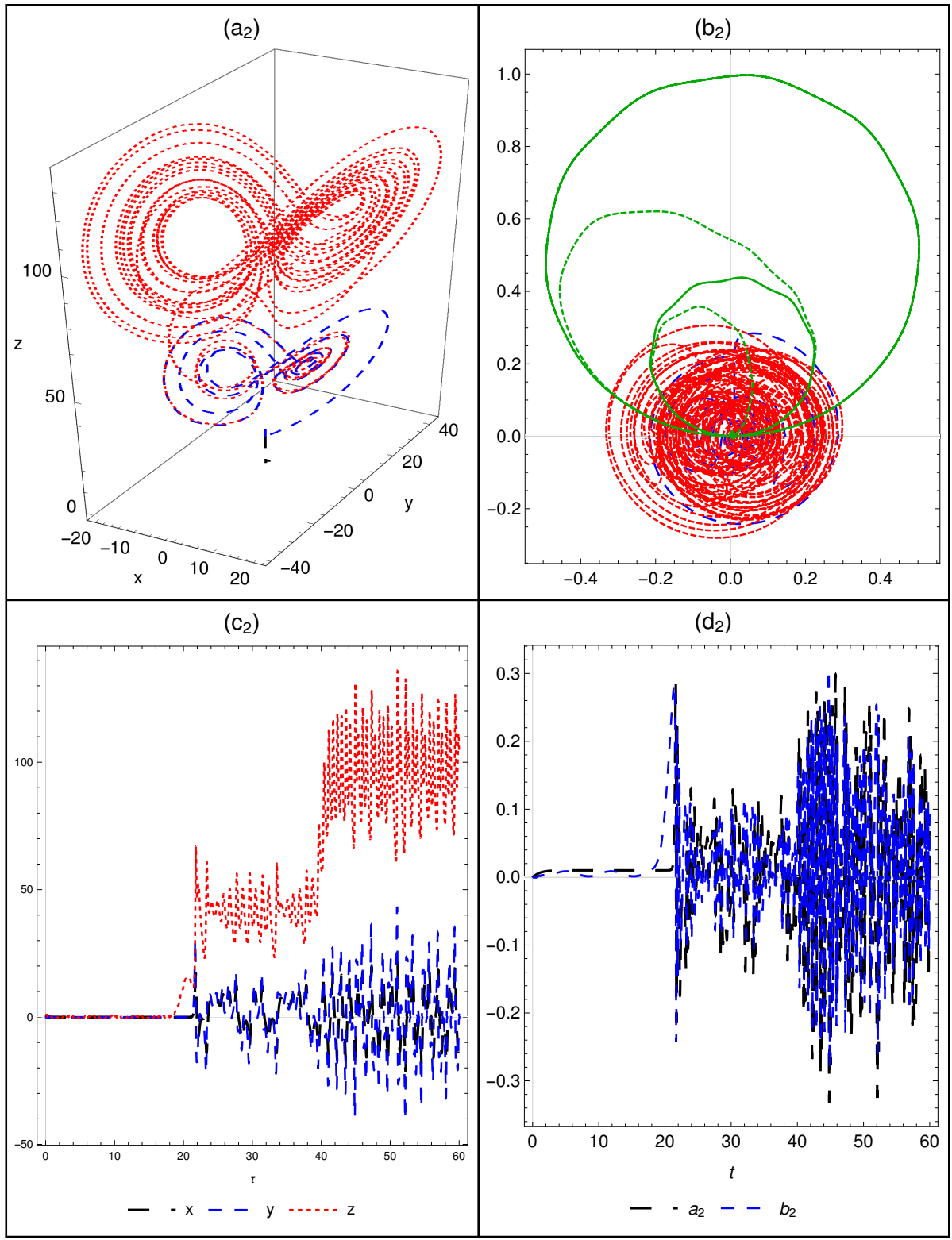}
\caption{Phase space and time series of solutions to \eqref{eq:aww-all-modes} for the unsteady symmetric water wheel when $x(0)=1$. In the top left panel we plot the phase portraits for $(x,y,z)\in\mathbb{R}^3$. In the top right panel we plot the phase space for the $n=2$ modes $a_2$ and $b_2$ as well as the evolving circle from \eqref{eq:evol-circle} (shown in green). We plot time series of $x(t)$, $y(t)$, and $z(t)$ in the lower left panel and of $a_2(t)$ and $b_2(t)$ in the lower right panel. }
 \label{fig:ps21}
\end{figure*}

\begin{figure*}
\centering
\includegraphics[width=0.9\textwidth]{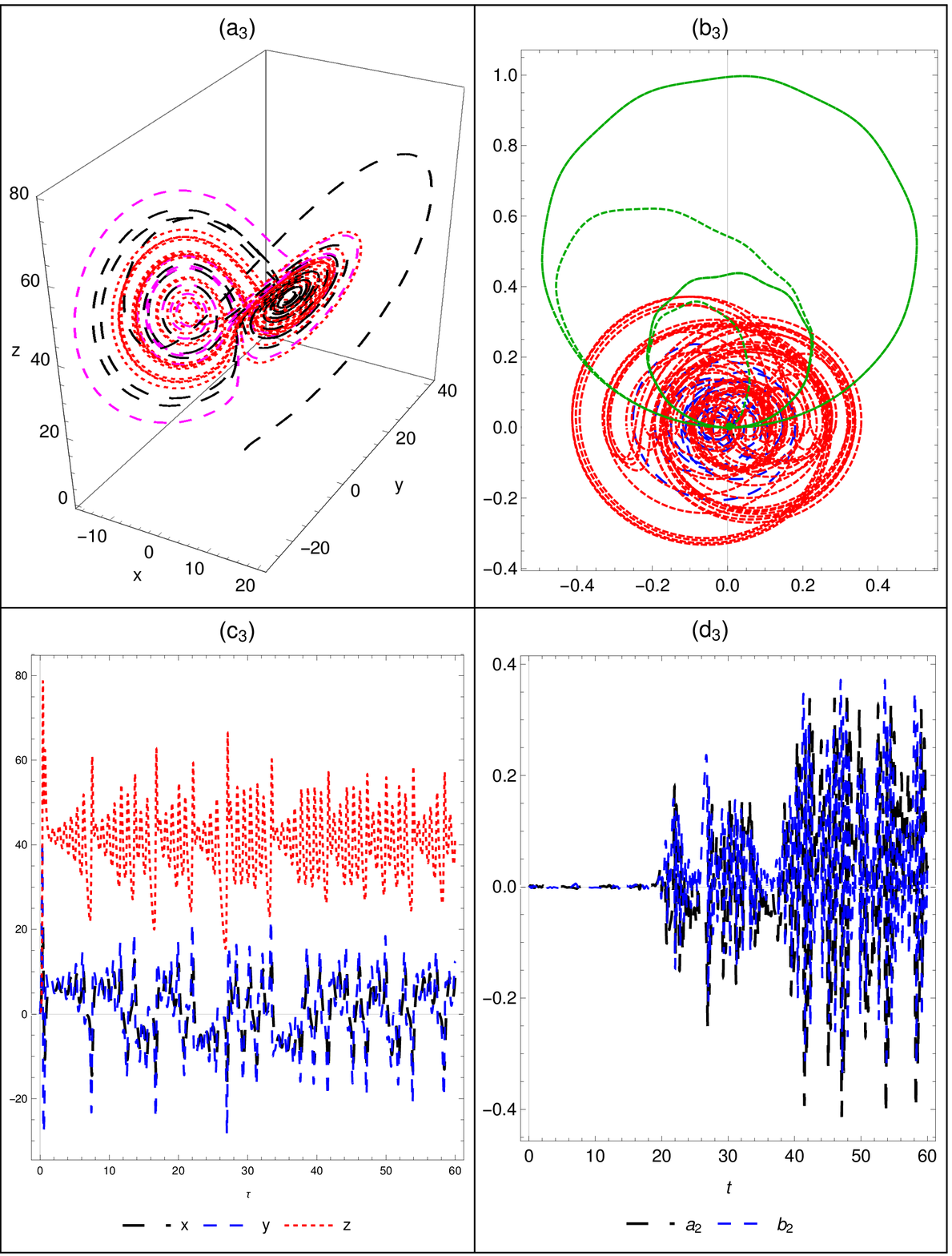}
\caption{Phase space and time series of solutions to \eqref{eq:aww-all-modes} for the steady asymmetric water wheel when $x(0)=1$. In the top left panel we plot the phase portraits for $(x,y,z)\in\mathbb{R}^3$. In the top right panel we plot the phase space for the $n=2$ modes $a_2$ and $b_2$ as well as the evolving circle from \eqref{eq:evol-circle} (shown in green). We plot time series of $x(t)$, $y(t)$, and $z(t)$ in the lower left panel and of $a_2(t)$ and $b_2(t)$ in the lower right panel. }
 \label{fig:ps31}
\end{figure*}

\begin{figure*}
\centering
\includegraphics[width=0.9\textwidth]{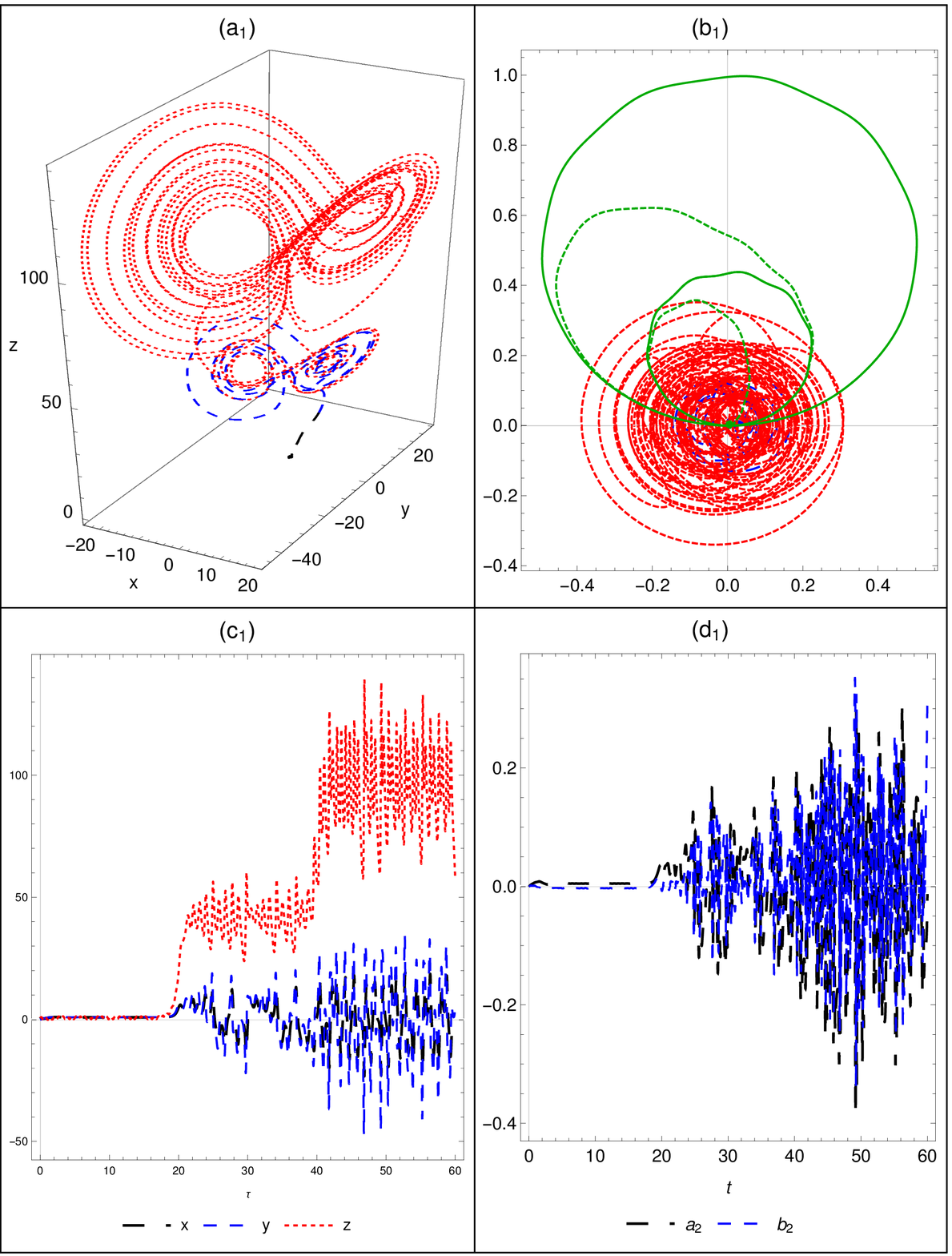}
\caption{Phase space and time series of solutions to \eqref{eq:aww-all-modes} for the unsteady asymmetric water wheel when $x(0)=1$. In the top left panel we plot the phase portraits for $(x,y,z)\in\mathbb{R}^3$. In the top right panel we plot the phase space for the $n=2$ modes $a_2$ and $b_2$ as well as the evolving circle from \eqref{eq:evol-circle} (shown in green). We plot time series of $x(t)$, $y(t)$, and $z(t)$ in the lower left panel and of $a_2(t)$ and $b_2(t)$ in the lower right panel. }
 \label{fig:ps11}
\end{figure*}

Figures \ref{fig:ps21}, \ref{fig:ps31}, and \ref{fig:ps11} show phase spaces and time series of the system variables obtained by solving \eqref{eq:taww} and \eqref{eq:aww-all-modes} numerically for the unsteady symmetric, steady asymmetric, and unsteady asymmetric water wheel models, respectively, taking $x(0)=1$. In the phase space plots, the black curve is the solution for $t\in [0,20]$, the blue curve is the solution for $t\in [20,40]$, and red curve is the solution for $t\in [40,60]$. In the top right panels, the green curve represents the evolving circle \eqref{eq:evol-circle}.

\begin{figure*}
\centering
\includegraphics[width=0.9\textwidth]{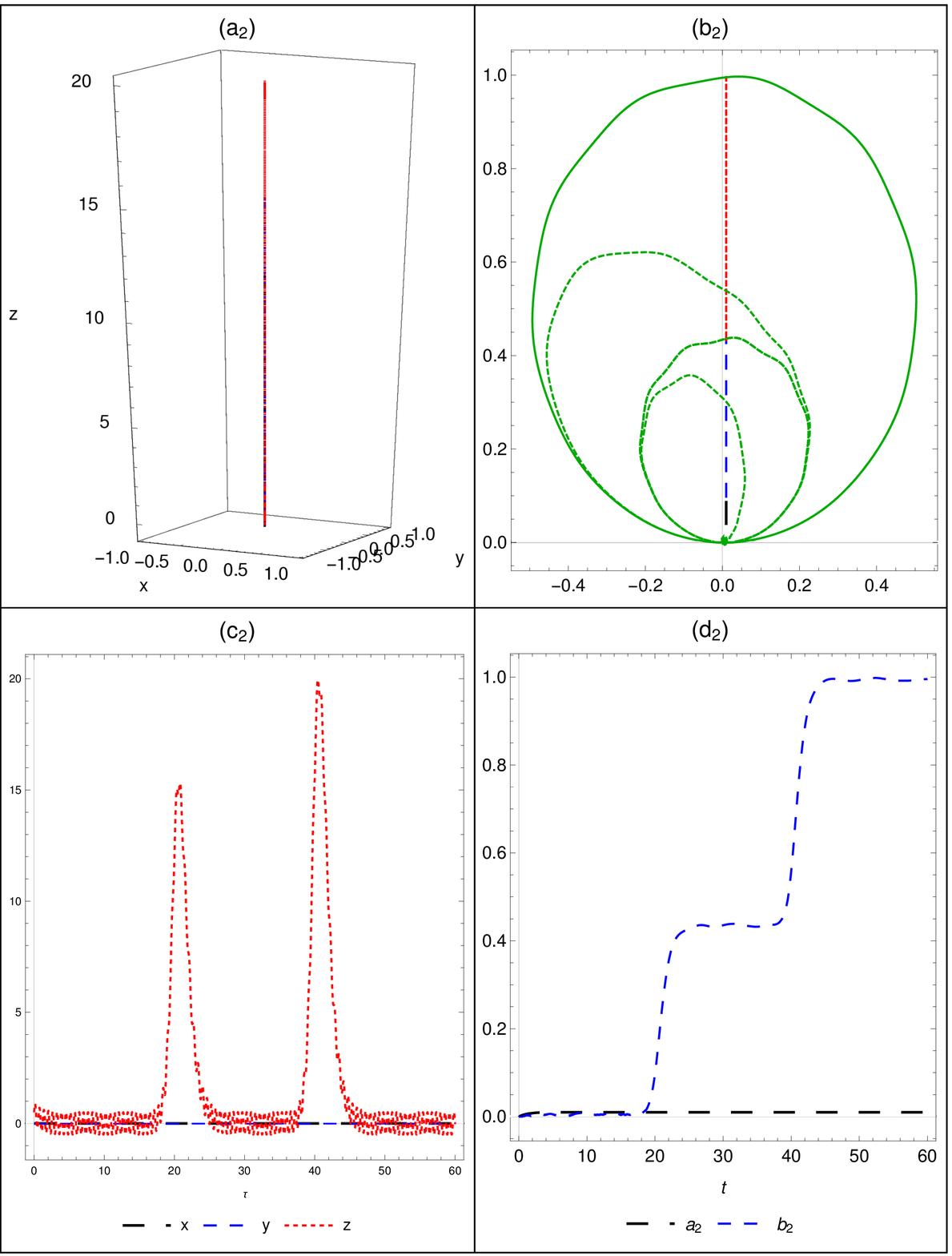}
\caption{Phase space and time series of solutions to \eqref{eq:aww-all-modes} for the unsteady symmetric water wheel when $x(0)=0$. In the top left panel we plot the phase portraits for $(x,y,z)\in\mathbb{R}^3$. In the top right panel we plot the phase space for the $n=2$ modes $a_2$ and $b_2$ as well as the evolving circle from \eqref{eq:evol-circle} (shown in green). We plot time series of $x(t)$, $y(t)$, and $z(t)$ in the lower left panel and of $a_2(t)$ and $b_2(t)$ in the lower right panel. }
 \label{fig:ps22}
\end{figure*}

\begin{figure*}
\centering
\includegraphics[width=0.9\textwidth]{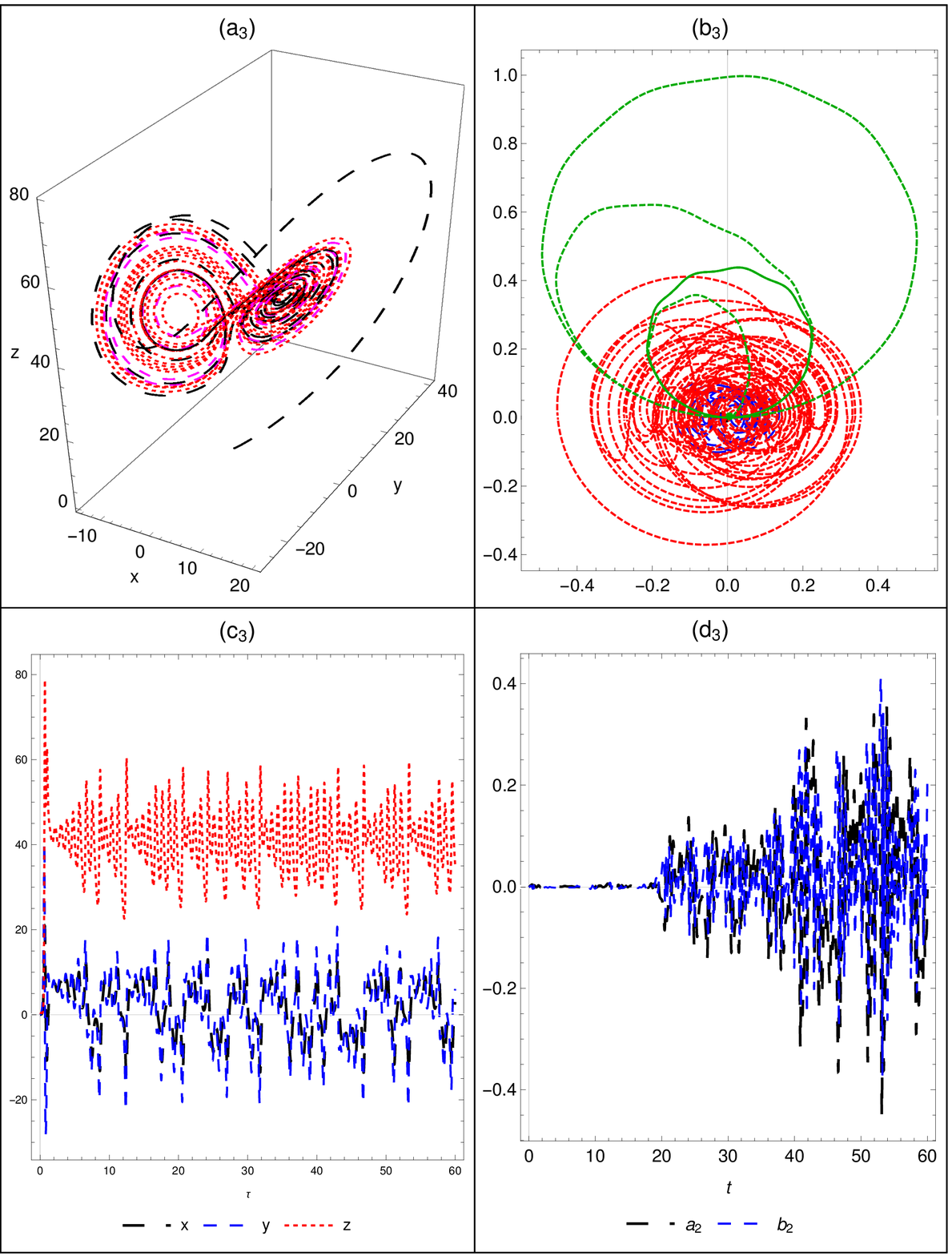}
\caption{Phase space and time series of solutions to \eqref{eq:aww-all-modes} for the steady asymmetric water wheel when $x(0)=0$. In the top left panel we plot the phase portraits for $(x,y,z)\in\mathbb{R}^3$. In the top right panel we plot the phase space for the $n=2$ modes $a_2$ and $b_2$ as well as the evolving circle from \eqref{eq:evol-circle} (shown in green). We plot time series of $x(t)$, $y(t)$, and $z(t)$ in the lower left panel and of $a_2(t)$ and $b_2(t)$ in the lower right panel.}
 \label{fig:ps32}
\end{figure*}

\begin{figure*}
\centering
\includegraphics[width=0.9\textwidth]{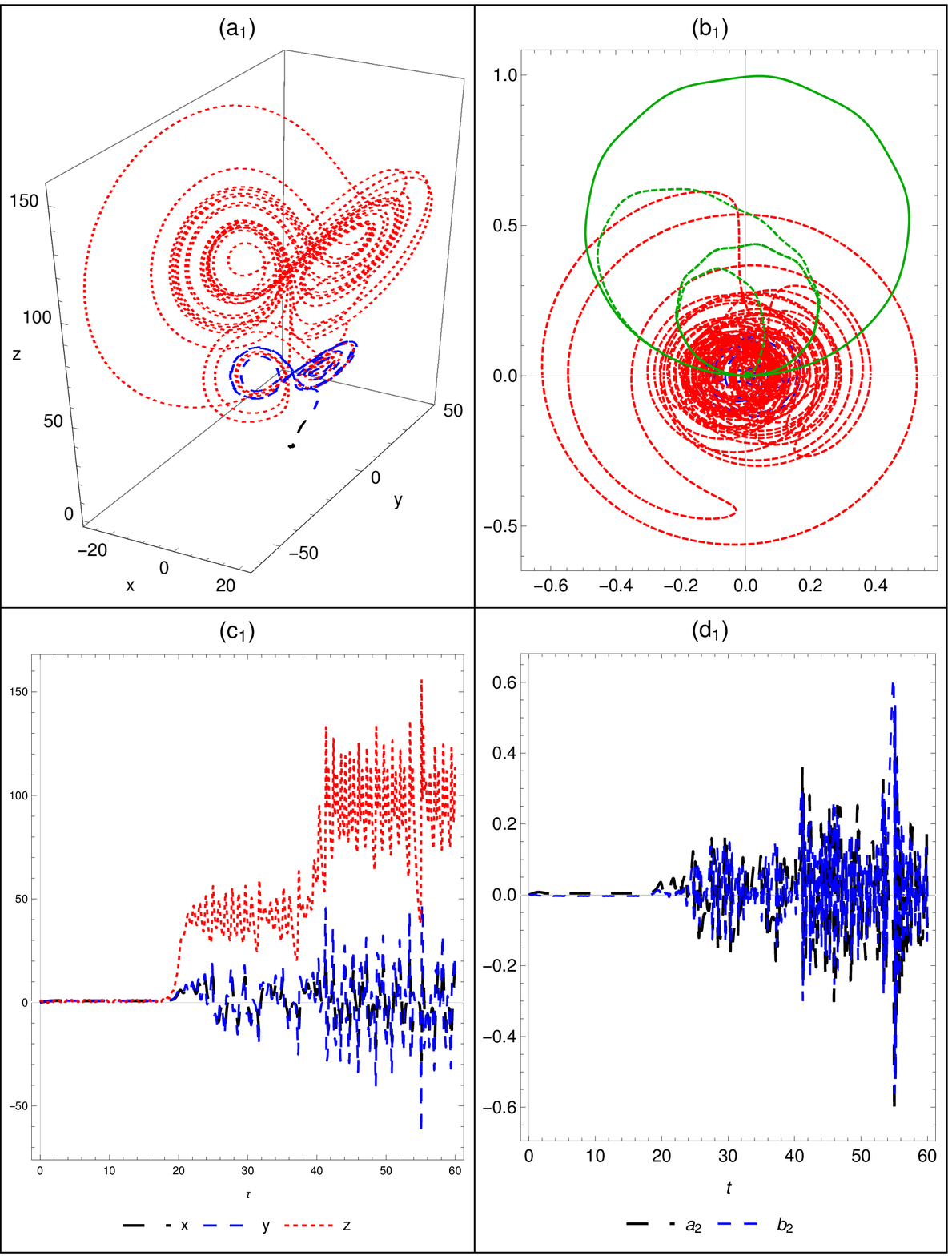}
\caption{Phase space and time series of solutions to \eqref{eq:aww-all-modes} for the unsteady asymmetric water wheel when $x(0)=0$. In the top left panel we plot the phase portraits for $(x,y,z)\in\mathbb{R}^3$. In the top right panel we plot the phase space for the $n=2$ modes $a_2$ and $b_2$ as well as the evolving circle from \eqref{eq:evol-circle} (shown in green). We plot time series of $x(t)$, $y(t)$, and $z(t)$ in the lower left panel and of $a_2(t)$ and $b_2(t)$ in the lower right panel. }
 \label{fig:ps12}
\end{figure*}

When the parameter $r(\tau)$ is non-zero, one expects chaos, as predicted by the competitive modes analysis. Since the parameters $r(\tau)$ and $\mu(\tau)$ are time dependent, the system adjusts its response and settles into chaotic motion in different phase space regions. For $t>20$, the $x$ and $y$ time-series switch signs unpredictably and the solution seems to settle on two different butterfly attractors. This indicates chaotic motion of the water wheel for large time. An examination of the phase portraits in Figures \ref{fig:ps21}-\ref{fig:ps11} and Figure \ref{fig:gijpos} reveals that a large part of the butterfly attractor lies outside of the red area in Figure \ref{fig:gijpos}. Also, the time series plot of the $z$ variable for all the water wheels satisfies equation of the manifold $g_2-g_3=0$ intermittently, thus satisfying condition 2 of the competitive modes Conjecture \ref{conj:chaos-cond}.

Figures \ref{fig:ps21}-\ref{fig:ps11} also show phase space trajectories of higher modes ($n=2$, although one can similarly plot modes for $n\geq 3$) in a chaotic regime and the evolving circle given in \eqref{eq:evol-circle} (the latter depicted by a green curve). As discussed in Section \ref{sec:qa}, all trajectories are attracted toward the boundary of the evolving circle \eqref{eq:evol-circle} for large time. Therefore, as the radius of the circle grows with time, so does the amplitude of the modes $a_n, b_n$. Also, since the phase space volume contracts (see \eqref{eq:vol-cont}), the modes remain bounded in a region of the phase space. 

We change the initial angular velocity to $x(0) =0$, and give corresponding results for the symmetric, steady asymmetric, and unsteady asymmetric water wheels in Figures \ref{fig:ps22}-\ref{fig:ps12}. For this modified condition, we do not expect the symmetric water wheel to spin at all (as discussed in Section \ref{sec:qa}), and this is the behavior which we observe in Figure \ref{fig:ps22}. On the other hand, the other two asymmetric water wheels can still move, and we observe again chaotic motion of these two water wheels in Figures \ref{fig:ps32} and \ref{fig:ps12}. In the absence of rotation in the symmetric water wheel case, modes corresponding to $n=2$ also get confined to a very small area of the phase space. This suggests that the modes occupy a larger proportion of perimeter of the circle in a chaotic regime compared to non-chaotic regimes.

\section{Conclusions} \label{sec:conc}
We have derived a non-autonomous nonlinear dynamical system governing an asymmetric water wheel with unsteady water inflow rate. In relevant limits, this reduces to a dynamical system for an asymmetric water wheel with steady water inflow rate, as well as a dynamical system for the symmetric water wheel. In order to determine parameter regimes likely to give chaotic dynamics in this model, we have applied the competitive modes analysis. As our nonlinear dynamical system is non-autonomous, we have modified the standard competitive modes approach by introducing an auxiliary function which accounts for the non-autonomous contribution, casting the system as a higher-order autonomous dynamical system for sake of the competitive modes analysis. Like the symmetric and asymmetric steady analogues, the asymmetric unsteady water wheel model permits chaotic dynamics for some parameter regimes. Owing to the non-autonomous contributions, we find that the dynamics are less regular in general, and an enlarged chaotic parameter regime is found when compared to the symmetric and asymmetric steady cases. Physically, our results suggest that chaos should be fairly ubiquitous in the asymmetric water wheel model with unsteady inflow of water, and can occur even in parameter regimes where the corresponding symmetric and asymmetric steady water wheel systems give regular, non-chaotic dynamics.

\end{document}